\documentclass[conference]{IEEEtran}

\usepackage{amsmath}
\usepackage{amsthm}
\usepackage{amssymb}
\usepackage{subcaption}
\usepackage{float}
\usepackage{turnstile}
\usepackage{mathtools}
\usepackage{cite}
\usepackage{color}
\usepackage{cite}
\usepackage{graphicx}

\newcommand\blfootnote[1]{%
  \begingroup
  \renewcommand\thefootnote{}\footnote{#1}%
  \addtocounter{footnote}{-1}%
  \endgroup
}

\newtheorem{theorem}{Theorem}

\newtheorem{assumption}{Assumption}

\newtheorem{lemma}{Lemma}

 \def\vr{{\bf r}}  
  \def\vw{{\bf w}} \def\vx{{\bf x}}

\def \bzeta{\boldsymbol{\zeta}}

\def \e{\varepsilon}
\def \R{\mathbb{R}}
\def \dx{\,dx}
\def \E{\mathbb{E}}

\def \N{\mathbb{N}}
\def \R{\mathbb{R}}
\def \calN{\mathcal{N}}

\mathtoolsset{showonlyrefs}

\begin{document}
\title{Distributed Global Optimization by Annealing}

\author{\IEEEauthorblockN{Brian Swenson\IEEEauthorrefmark{1},
Soummya Kar\IEEEauthorrefmark{2},
H. Vincent Poor\IEEEauthorrefmark{1}, and
Jos{\'e} M. F. Moura\IEEEauthorrefmark{2}}
\IEEEauthorblockA{\IEEEauthorrefmark{1}Department of Electrical Engineering\\
Princeton University,
Princeton, NJ 08540\\ Email: bswenson@princeton.edu, poor@princeton.edu}
\IEEEauthorblockA{\IEEEauthorrefmark{2}Department of Electrical and Computer Engineering\\
Pittsburgh, PA 15213\\
Email: soummyak@andrew.cmu.edu, moura@andrew.cmu.edu}}


\maketitle

\begin{abstract}
The paper considers a distributed algorithm for global minimization of a nonconvex function. The algorithm is a first-order \emph{consensus + innovations} type algorithm that incorporates decaying additive Gaussian noise for annealing, converging to the set of global minima under certain technical assumptions. The paper presents simple methods for verifying that the required technical assumptions hold and illustrates it with a distributed target-localization application.
\end{abstract}


\thispagestyle{empty}
\markboth{}{}

\section{Introduction}
\label{introduction}
\blfootnote{The work of B. Swenson and H. V. Poor was partially supported by the Air Force Office of Scientific Research under MURI Grant FA9550-18-1-0502.
The work of S. Kar and J. M. F. Moura was partially supported by the National Science Foundation (NSF) under NSF Grant CCF 1513936.}
Nonconvex optimization problems are prevalent throughout machine learning and signal processing \cite{goodfellow2016deep,chartrand2007exact,chen2018gradient,kar2019clustering}.
In settings, such as the internet of things (IoT) and sensor networks, it can be impractical to process information in a centralized fashion due to the high volume of data inherently distributed across many devices \cite{shi2016edge,hartenstein2008tutorial,hu2015mobile}.
Moreover, due to privacy concerns, users can be unwilling to share (potentially sensitive) data for processing in a central location \cite{abadi2016deep}.
This necessitates the development of \emph{distributed} algorithms for nonconvex optimization.

We are interested in studying distributed algorithms wherein (i) agents may only communicate with neighbors via an overlaid communication network (possibly time-varying), and (ii) there is no central node or entity to coordinate the computation. Within this framework, we consider distributed algorithms to optimize the function
\begin{equation} \label{eq:U-def}
U(x) := \sum_{n=1}^N U_n(x),
\end{equation}
where $N$ denotes the number of agents in the network and $U_n$ is a local function available only to agent $n$. Example applications of distributed nonconvex optimization problems in this framework include empirical risk minimization \cite{lee2018distributed}, target localization \cite{di2015distributed}, robust regression \cite{sun2016distributed}, distributed coverage control \cite{welikala2019distributed}, power allocation in wireless adhoc networks \cite{bianchi2012convergence}, and others \cite{di2016next}.

Work on distributed nonconvex optimization has focused largely on ensuring convergence to first-order stationary points \cite{bianchi2012convergence,zhu2012approximate,magnusson2015convergence,di2016next,
sun2016distributed,tatarenko2017non}.
More recently, \cite{daneshmand2018second,daneshmand2018second-b,hong2018gradient} have considered the problem of demonstrating convergence to local optima and evasion of saddle points. In this paper, we consider the problem of developing distributed algorithms for computing global optima.

We will focus on the following annealing-based algorithm:
\begin{align} \label{eq:update}
\mathbf{x}_{n}(t+1) =\, & \mathbf{x}_{n}(t)-\beta_{t}\sum_{\ell\in\Omega_{n}(t)}\left(\mathbf{x}_{n}(t)-\mathbf{x}_{\ell}(t)\right)\\
& -\alpha_{t}\left(\nabla U_{n}(x_{n}(t))+\bzeta_{n}(t)\right)+\gamma_{t}\mathbf{w}_{n}(t),
\end{align}
$n=1,\ldots,N$, where $\vx_n(t)\in\R^d$ is the state of agent $n$ at iteration $t\geq 1$; $\Omega_n(t)$ denotes the set of agents neighboring $n$ at time $t$ (per the communication graph); $\{\alpha_t\}$, $\{\beta_t\}$, and $\{\gamma_t\}$ are sequences of decaying weight parameters; $\zeta_n$(t) is a $d$-dimensional random variable (representing gradient noise); and $\vw_n(t)$ is $d$-dimensional Gaussian noise (introduced for annealing).
Note that this algorithm is distributed in the sense that, to compute $\vx_n(t+1)$, each agent only requires information about their local function $U_n$ and the state $\vx_\ell(t)$ of neighboring agents $\ell\in \Omega_n(t)$.

Algorithm \eqref{eq:update} is a \emph{consensus + innovations} type algorithm \cite{kar2012distributed}. The first term $-\beta_{t}\sum_{l\in\Omega_{n}(t)}\left(\mathbf{x}_{n}(t)-\mathbf{x}_{l}(t)\right)$
(referred to as the \emph{consensus term}) ensures that agents reach asymptotic agreement; the second term $-\alpha_{t}\left(\nabla U_{n}(x_{n}(t))+\bzeta_{n}(t)\right)$ (referred to as the \emph{innovations term}) ensures that agents descend their local objective $U_n$; the final term $\gamma_{t}\mathbf{w}_{n}(t)$ is an annealing term that ensures that limit points are global rather than local minima.
The algorithm may be seen as a distributed variant of the (centralized) annealing-based algorithm studied in \cite{gelfand1991recursive}.

Convergence properties of \eqref{eq:update} were studied in \cite{swenson2019CDC} where, under certain assumptions, it was established that the algorithm converges in distribution to the set of global optima of \eqref{eq:U-def}.
Some of the assumptions under which this convergence result is proved are highly technical.  In this paper we will review the convergence results for \eqref{eq:update} and present simple methods for verifying that the required assumptions hold in the context of a target localization example.

The remainder of the paper is organized as follows. Section~\ref{sec:notation} sets up notation. Section~\ref{sec:conv-result} presents our assumptions and the convergence result for \eqref{eq:update}. Section~\ref{sec:example} discusses a distributed target localization application. Section~\ref{sec:conclusion} concludes the paper.

\section{Notation} \label{sec:notation}
\label{notgraph}
We use $\|\cdot\|$ to indicate the standard Euclidean norm. Given $x\in \R^d$ and $r>0$, $B_r(x)$ denotes the open ball of radius $r$ about $x$. We let $\lambda$ denote the Lebesgue measure and let $\N:=\{1,2,\ldots\}$. For $k\geq 1$, we say a function $f:\R^d\to\R$ is of class $C^k$ if $f$ is $k$-times continuously differentiable. For a function $f:\R^d\to\R$, when well defined, we let $\nabla f(x)$ denote the gradient, $\nabla^2 f(x)$ denote the Hessian, and $\Delta f(x) = \sum_{i=1}^d \frac{\partial^2 f(x)}{\partial x_i^2}$ denote the Laplacian of $f$.

We will assume that agents may communicate over an undirected, time-varying graph $G_t = (V_t,E_t)$, where $V_t$ denotes the set of vertices (or agents) and $E_t$ denotes the set of edges. We assume that $G_t$ is devoid of self-loops, so that $(i,i)\notin E_t$ for any $i\in\{1,\ldots,N\}$ or $t\geq 1$. A link $(i,j)\in E_t$ denotes the ability of agents $i$ and $j$ to communicate at time $t$. The adjacency matrix associated with $G_t$ is given by $A = (a_{ij}^t)$, where $a_{ij}^t = 1$ if $(i,j)\in E_t$, and $a_{ij}^t=0$ otherwise, and the degree matrix is given by the diagonal matrix $D_t$ with diagonal entries $d_i^t = \sum_{j=1}^N a_{ij}^t$. The graph Laplacian of $G_t$ is given by the matrix $L_t = D_t - A_t$.

Given a measure $\pi$ on a measurable space $(\mathbb{R}^{k},\Sigma)$, and a (measurable) function $f:\mathbb{R}^{k}\mapsto\mathbb{R}$, we use the convention
\begin{equation}\label{eq:pi-f-def}
\pi(f) ~ := \int f d\pi,
\end{equation}
whenever the integral exists.
For a stochastic process $\{Z_{t}\}$, with $Z_t\in \R^d$, and a function $f:\R^d\to\R$, we use the convention
\begin{equation} \label{eq:def-cond-E}
\E_{t_0,z_0}[f(Z_t)] ~ := ~ \E[f(Z_t)\vert Z_{t_0} = z_0].
\end{equation}

\section{Convergence Result} \label{sec:conv-result}
\subsection{Assumptions} \label{sec:ass}
The main convergence result for \eqref{eq:update} will be given in Theorem~\ref{th:main_result}. We will make the following assumptions.
\begin{assumption}\label{ass:Lip}
The function $U_{n}(\cdot)$ is $C^{2}$ with Lipschitz continuous gradients, i.e., there exists an $L>0$ such that
\begin{align}
\left\|\nabla U_{n}(x)-\nabla U_{n}(\acute{x})\right\|\leq L\left\|x-\acute{x}\right\|,  \quad \forall n.
\end{align}
\end{assumption}
\begin{assumption}
\label{ass:similarity} The function $U_{n}(\cdot)$ satisfies the following bounded gradient-dissimilarity condition:
\begin{align}
\sup_{x\in\mathbb{R}^{d}}\left\|\nabla U_{n}(x)-\nabla U(x)\right\|<\infty, \quad \forall n.
\end{align}
\end{assumption}

\begin{assumption}
\label{ass:GM_1} $U(\cdot)$ is a $C^2$ function such that
\begin{enumerate}
\item[(i)] $\min U(x) = 0$,
\item[(ii)] $U(x)\to\infty$ and $\|\nabla U(x)\|\to\infty$ as $\|x\|\to\infty$,
\item[(iii)] $\inf(\|\nabla U(x)\|^2 - \Delta U(x) ) > -\infty$.
\end{enumerate}
\end{assumption}
\begin{assumption}
\label{ass:GM_2} For $\e>0$ let
\begin{align}
d\pi^\e(x) & = \frac{1}{Z^\e}\exp\left(-\frac{2U(x)}{\e^2} \right)\dx,\\
Z^\e & = \int\exp\left(-\frac{2U(x)}{\e^2} \right)\dx.
\end{align}
$U$ is such that $\pi^\e$ has a weak limit $\pi$ as $\e\to 0$.
\end{assumption}
We note that $\pi$ is constructed so as to place mass 1 on the set of global minima of $U$. A simple condition ensuring the existence of such a $\pi$ will be discussed in Lemma~\ref{lemma:Laplace-ref}.
\begin{assumption}
\label{ass:GM_3} $\liminf_{\|x\|\to\infty}\left\langle \frac{\nabla U(x)}{\|\nabla U(x)\|}, \frac{x}{\|x\|} \right\rangle \geq C(d)$,
$C(d) = \left( \frac{4d-4}{4d-3} \right)^{\frac{1}{2}}$.
\end{assumption}
\begin{assumption}
\label{ass:GM_5} $\liminf_{\|x\|\to\infty} \frac{\|\nabla U(x)\|}{|x|} > 0$
\end{assumption}
\begin{assumption}
\label{ass:GM_6}
$\limsup_{\|x\|\to\infty} \frac{\|\nabla U(x)\|}{\|x\|} < \infty$
\end{assumption}

Let $\{\mathcal{H}_{t}\}$ denote the natural filtration corresponding to the update process~\eqref{eq:update}, i.e., for all $t$,
\begin{align}
\label{eq:filt}
\mathcal{H}_{t}=\sigma\left(\mathbf{x}_{0},L_{0},\cdots,L_{t-1},\bzeta_{0},\cdots,\bzeta_{t-1},\mathbf{w}_{0},\cdots,\mathbf{w}_{t-1}\right).
\end{align}

\begin{assumption}
\label{ass:conn} The $\{\mathcal{H}_{t+1}\}$-adapted sequence of undirected graph Laplacians $\{L_{t}\}$ are independent and identically distributed (i.i.d.), with $L_{t}$ being independent of $\mathcal{H}_{t}$ for each $t$, and are connected on the mean, i.e., $\lambda_{2}(\bar{L}_t)>0$ where $\bar{L}_t=\mathbb{E}[L_{t}]$.
\end{assumption}

\begin{assumption}
\label{ass:grad_noise} The sequence $\{\bzeta_{t}\}$ is $\{\mathcal{H}_{t+1}\}$-adapted and there exists a constant $L_{1}>0$ such that
\begin{align}
\label{ass:grad_noise1}
\mathbb{E}[\bzeta_{t}~|~\mathcal{H}_{t}]=0~~\mbox{and}~~\mathbb{E}[\|\bzeta_{t}\|^{2}~|~\mathcal{H}_{t}]<L_{1}
\end{align}
for all $t\geq 0$.
\end{assumption}
\begin{assumption}
\label{ass:gauss}
For each $n$, the sequence $\{\mathbf{w}_{n}(t)\}$ is a sequence of i.i.d. $d$-dimensional standard Gaussian vectors with covariance $I_d$ and with $\mathbf{w}_{n}(t)$  being independent of $\mathcal{H}_{t}$ for all $t$. Further, the sequences $\{\mathbf{w}_{n}(t)\}$ and $\{\mathbf{w}_{\ell}(t)\}$ are mutually independent for each pair $(n,\ell)$ with $n\neq \ell$.
\end{assumption}
\begin{assumption}
\label{ass:weights} The sequences $\{\alpha_{t}\}$, $\{\beta_{t}\}$, and $\{\gamma_{t}\}$ satisfy
\begin{align}
\label{eq:weights}
\alpha_{t}=\frac{c_{\alpha}}{t},~~\beta_{t}=\frac{c_{\beta}}{t^{\tau_{\beta}}},~~\gamma_{t}=\frac{c_{\gamma}}{t^{1/2}\sqrt{\log\log t}},~~~\mbox{for $t$ large},
\end{align}
where $c_{\alpha},c_{\beta},c_{\gamma}>0$ and $\tau_{\beta}\in (0,1/2)$.
\end{assumption}

Finally, let $C_{0}$ be the constant as defined after (2.3) in~\cite{gelfand1991recursive}.

Assumptions~\ref{ass:Lip}--\ref{ass:similarity} ensure that each $U_n$ is well behaved so that asymptotic consensus can be achieved by the consensus component of \eqref{eq:update}. Assumptions \ref{ass:GM_1}--\ref{ass:GM_6} ensure that $U$ is well behaved so that convergence to a global minimum is possible. Assumption \ref{ass:conn} ensures that the communication network is sufficiently well connected, while Assumptions \ref{ass:grad_noise} and \ref{ass:gauss} ensure that the algorithm explores the state space and tends towards a descent direction. Finally, Assumption \ref{ass:weights} ensures that the weight parameters in \eqref{eq:update} appropriately balance the objectives of reaching consensus, descending the objective, and exploring the state space.

\subsection{Convergence Result}
The main convergence result for process \eqref{eq:update} is given in the following theorem. Informally, the theorem states that $X_n$ converges in distribution to a random variable placing mass 1 on the set of global minima of $U$.

\begin{theorem}
\label{th:main_result}
Let $\{\vx_1(t)\}_{t\geq 1},\ldots,\{\vx_N(t)\}_{t\geq 1}$ satisfy the recursion \eqref{eq:update} with respective initial conditions $x_{0,1},\ldots,x_{0,N}$. Let Assumptions~\ref{ass:Lip}--\ref{ass:weights} hold.
Further, suppose that $c_{\alpha}$ and $c_{\gamma}$ in Assumption~\ref{ass:weights} satisfy $c_{\gamma}^{2}/c_{\alpha}>C_{0}$, where $C_0$ is defined after Assumption~\ref{ass:weights}.
Then, for any bounded continuous function $f:\mathbb{R}^{d}\to\mathbb{R}$ and for all $n=1,\ldots,N$, we have that
\begin{align}
\label{eq:main_result1}
\lim_{t\rightarrow\infty}\mathbb{E}\left[f(\mathbf{x}_{n}(t))\right]=\pi(f).
\end{align}
\vspace{-2em}
\end{theorem}
In the above theorem we recall that we use conventions \eqref{eq:pi-f-def}--\eqref{eq:def-cond-E}. We also recall that the relationship between the condition for weak convergence used in \eqref{eq:main_result1} and other conditions for weak convergence (or convergence in distribution) are elucidated in the so-called \emph{portmanteu theorem} \cite{billingsley2013convergence}. A complete proof of Theorem~\ref{th:main_result} can be found in \cite{swenson2019CDC}.

\section{Illustrative Example} \label{sec:example}
Consider the problem of using $N$ sensors in a network to collaboratively locate the position of $T$ targets lying on a plane. All sensors and targets lie within some compact set $K$ with diameter $R=\max_{x_1,x_2\in K}\|x_1-x_2\|$, known apriori to all agents. Each sensor $n$ has knowledge of its own location $s_n\in \R^2$ and the distance between itself and target $k$, denoted by $d_{nk}$.

In order to formulate the problem of localizing the targets as an optimization problem, we begin by defining the following auxiliary function.
Given arbitrary $y,r\in [0,\infty)$, let
$$
g(y,r) :=
\begin{cases}
\phi_1(y,r) & 0\leq y \leq r/2\\
(y-r)^2 & r/2 < y
\end{cases}
$$
where, for all $r\in[0,\infty)$, the function $y\mapsto \phi_1(y,r)$ is finite-valued, monotone increasing, of class $C^3$, and is chosen such that $y\mapsto g(y,r)$ is also $C^3$ ($\phi_1$ may be constructed using a Hermite interpolating polynomial \cite{stoer2013introduction}). Given $x\in\R^2$ and $d\in [0,\infty)$ let
$$
f(x,r) := g(\|x\|,r).
$$
\vspace{-1.5em}

\begin{figure}[h]
  \begin{subfigure}[h]{0.20\textwidth}
    \includegraphics[width=\textwidth]{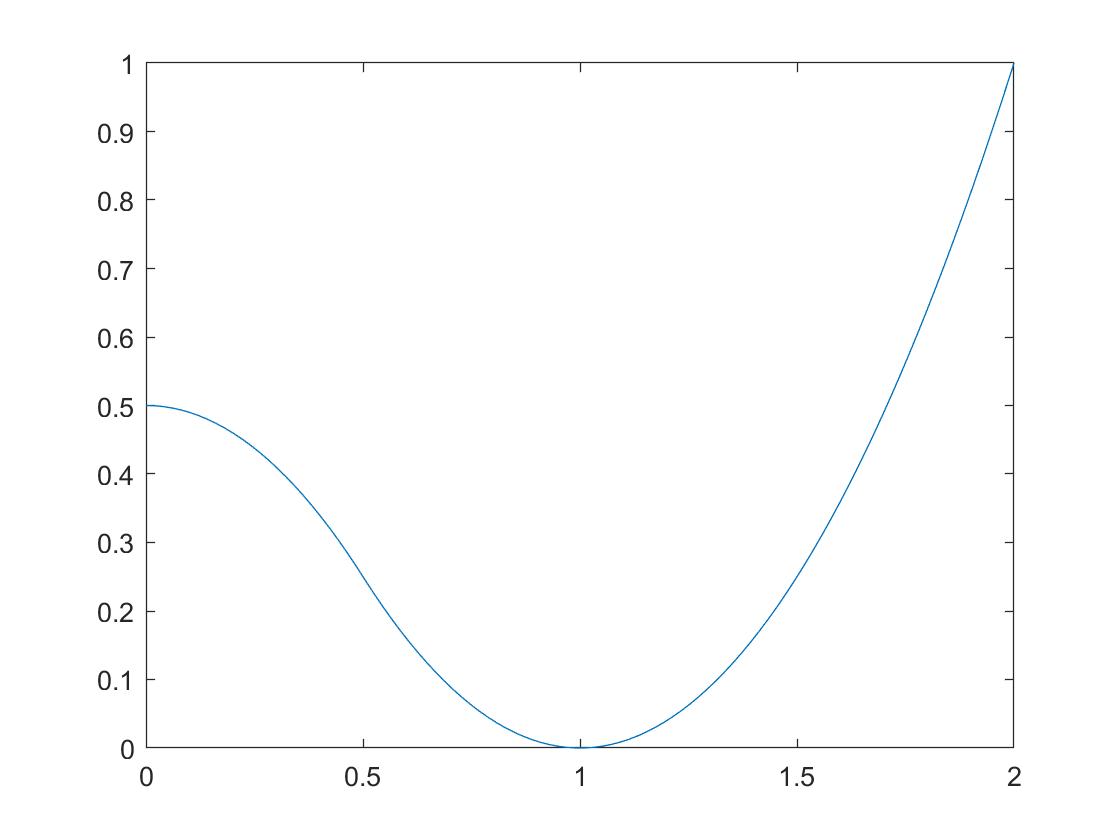}
    \caption{Example of function $f$ used in the construction of $g$.}
    \label{fig:g_fun}
  \end{subfigure}
  $~$
  \begin{subfigure}[h]{0.25\textwidth}
    \includegraphics[width=\textwidth]{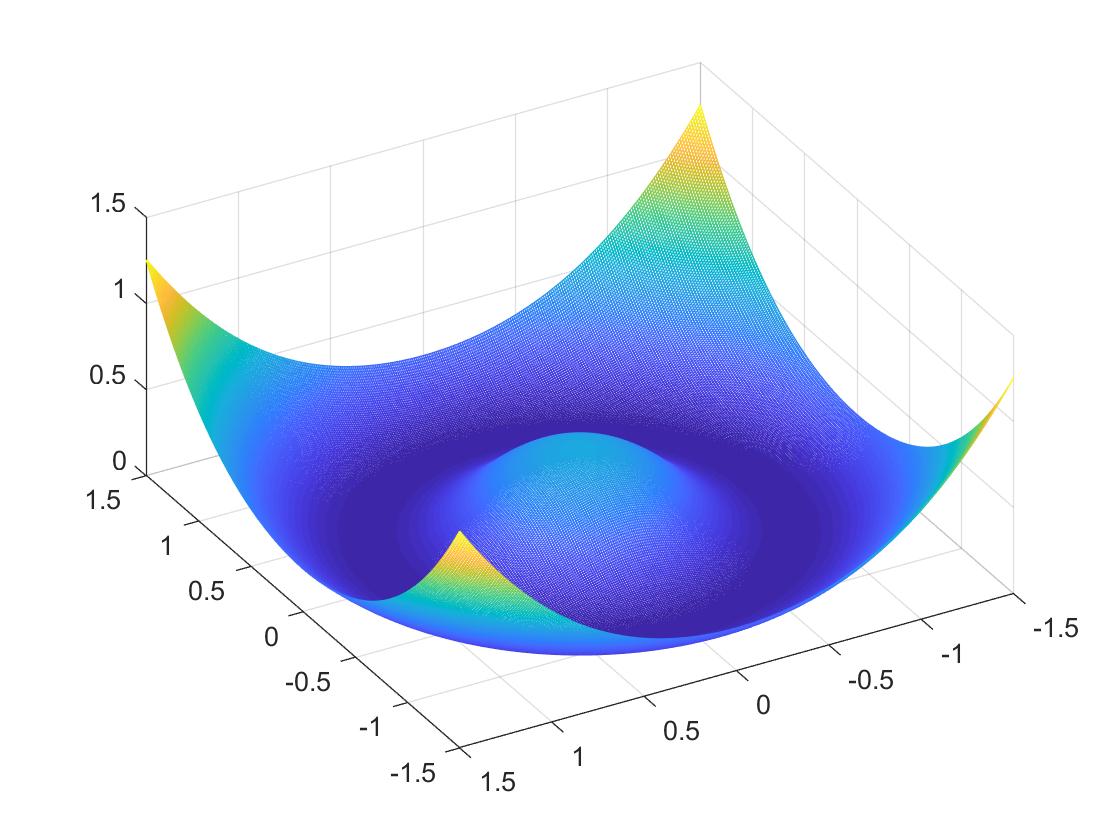}
    \caption{Example of function $g$ used in the construction of $U_n$.}
    \label{fig:f_fun}
  \end{subfigure}
\end{figure}

Examples of the functions $g$ and $f$ are shown in Figures~\ref{fig:g_fun} and \ref{fig:f_fun}. These functions are prototypes that will be used in the formulation of the optimization problem. Informally, if $s$ is the location of a sensor and it is known that a target lies a distance $r$ from $s$, then $f(x-s,r)$ is minimized (with value zero) along the ring with radius $r$ about $s$. Finally, we let the objective of player $n$ be given by
\begin{equation} \label{def:U_n-ex}
U_n(x) =
\begin{cases}
\sum_{k=1}^T f(x_n^k - s_n, d_{nk}) & \frac{1}{T}\|x_n\|< R\\
\phi_2(x_n) & R\leq \frac{1}{T}\|x_n\| < R+1\\
\|x_n\|^2 & \frac{1}{T}\|x_n\| > R,
\end{cases}
\end{equation}
where $x^k\in \R^2$ is an estimate of the location of target $k$, $x$ is the vector stacking $(x^k)_{k=1}^T$, and $\phi_1$ and $\phi_2$ are chosen so that
$U_n$ is $C^3$. The functions $\phi_1$ and $\phi_2$ may be constructed using Hermite interpolating polynomials \cite{stoer2013introduction}.

The function $U_n$ may be interpreted as follows. For $x$ with $x^k\in K$, $k=1,\ldots,T$ the function $U_n$ operates ``as expected'' (assigning high cost if $\|x^k - s_n\|$ is not close to $d_{nk}$).
For $x^k$ outside the set $K$, we have $U_n(x) = \|x\|^2$ for all $n$. This ensures that Assumption~\ref{ass:similarity} is satisfied. Finally, the transitory component $\phi_2$ merely ensures that $U_n$ transitions sufficiently smoothly between these two modes.\footnote{We remark that similar formulations of this problem using a quartic objective function have been considered in \cite{chen2012diffusion,di2016next}. Here, we reformulate the problem in terms of a quadratic objective in order to ensure that the assumptions in Section~\ref{sec:ass} are satisfied. See Section~\ref{sec:verification} for more details.}

Given $U_n$, $n=1,\ldots,N$ as defined in \eqref{def:U_n-ex}, the target localization problem is formulated as the unconstrained optimization of the sum function \eqref{eq:U-def}.

\subsection{Verification of Assumptions} \label{sec:verification}
In Theorem~\ref{th:main_result} we assumed that $U$ and $U_n$ satisfied Assumptions~\ref{ass:Lip}--\ref{ass:GM_6}.\footnote{The remaining Assumptions~\ref{ass:conn}--\ref{ass:weights} concern the algorithm \eqref{eq:update}, and not the optimization problem.} We now verify that these assumptions hold in the target localization example.

Assumptions~\ref{ass:Lip}--\ref{ass:similarity} hold since, by construction, $U_n(x) = \|x\|^2$ for $x$ sufficiently large. It is straightforward to verify that parts (i)--(ii) of Assumption~\ref{ass:GM_1} hold. Part (iii) of Assumption~\ref{ass:GM_1} holds due to the fact that each $U_n$ is quadratic for $x$ sufficiently large. In particular, for $x$ large we have $U_n(x) = \|x\|^2$, so that $\|\nabla U_n(x)\|^2 \geq \Delta U_n(x)$.

The following result from \cite{hwang1980laplace} (see \cite{hwang1980laplace}, Theorem 3.1) will allow us to verify that Assumption~\ref{ass:GM_2} holds.\footnote{Lemma~\ref{lemma:Laplace-ref} below has been adapted to fit our presentation and is slightly weaker than the result proved in \cite{hwang1980laplace}.}
\begin{lemma} \label{lemma:Laplace-ref}
Let
$\calN := \{x:U(x) = \inf_x U(x)\}.$
Suppose that \\
\noindent (i) $\lambda(\{U(x)< a \}) > 0$ for any $a> \inf_x U(x)$,\\
\noindent (ii) $\min_x U(x)$ exists and equals zero,\\
\noindent (iii) There exists $\e>0$ such that $\{U(x) \leq \e\}$ is compact,\\
\noindent (iv) $U$ is $C^3$.\\
Assume that $\calN$ consists of a finite set of isolated points and that the Hessian $\nabla^2 U(x)$ is invertible for all $x\in \calN$. Then the limit $\pi$ in Assumption~\ref{ass:GM_3} exists.
\end{lemma}

It is straightforward to verify that the sum function $U$ satisfies conditions (i)--(iv) of Lemma~\ref{lemma:Laplace-ref}. To verify that the Hessian $\nabla^2 U(x)$ is invertible for $x\in \calN$ requires some additional care.

Letting, $\{z_k\}_{k=1}^T$, $z_k\in \R^2$ denote the set of target locations, we will make the following additional assumption.
\begin{assumption}\label{ass:example-degen}$~$\\
\noindent (i) $N\geq 3$ and $s_n \not= s_m$ for all sensors $n\not =m$. \\
\noindent (ii) For each target $k$, there exist at least two sensors $m$ and $n$ such that $z_k$, $s_n$, and $s_m$ are not colinear. That is, $z_k \not = s_n + \lambda(s_m-s_n)$ for any $\lambda\in \R$.
\end{assumption}
Under part (i) of this assumption, the vector of targets $z=(z_k)_{k=1}^T$ is the unique global minimum of $U$, i.e., $\calN=\{z\}$. Under part (ii) of this assumption, the Hessian $\nabla^2 U(x)$ is invertible at $z$. This may be confirmed algebraically using the form of $U_n$ in \eqref{def:U_n-ex}.

Finally, Assumptions~\ref{ass:GM_3}--\ref{ass:GM_6} are seen to hold by again using the fact that each $U_n$ is quadratic for $x$ large.

In the numerical example to be given next, we will explicitly choose the graph $G_t$, weight sequences $\{\alpha_t\}$, $\{\beta_t\}$, and $\{\gamma_t\}$, and random variables $\bzeta_n(t)$ and $\vw_n(t)$ so that the remaining assumptions (Assumptions~\ref{ass:conn}--\ref{ass:weights}) are satisfied.

\subsection{Numerical Example}
In this section we consider a simple numerical example illustrating the functioning of the distributed annealing algorithm. We emphasize that these results are not optimal---the parameters are not chosen to optimize convergence rate, but merely to illustrate the general operation of the algorithm.

Consider an example of the target localization problem having five sensors and one target.\footnote{The choice of small parameter sizes for $N$ and $T$ in this example facilitates the visualization of the algorithm by allowing us to visually relate the behavior to the asymptotic mean vector field with global and local minima in simple figures.} The sensors are connected via a ring graph. The function $\phi_1$ (used in $g$) is constructed using a Hermite polynomial to smoothly interpolate between the functions $(x-d)^2$ (outside the ($d/2$)-ball) and $-x^2+(d/2)^2$ (inside a ($d/2-\e$)-ball).\footnote{In these simulations, we used $U_n(x) = \sum_{k=1}^T f\left(\|x_n^k - s_n\| - d_{nk}\right)$. In all simulations, trajectories $x_n^k(t)$ remained in the ball $B_3(0)$, so incorporating the remaining components of $U_n$ in \eqref{def:U_n-ex} was unnecessary. An interesting future research direction may be to formally relax Assumption~\ref{ass:similarity} in Theorem~\ref{th:main_result}.}

\begin{figure}[h]
  \begin{subfigure}[h]{0.45\textwidth}
    \includegraphics[width=\textwidth]{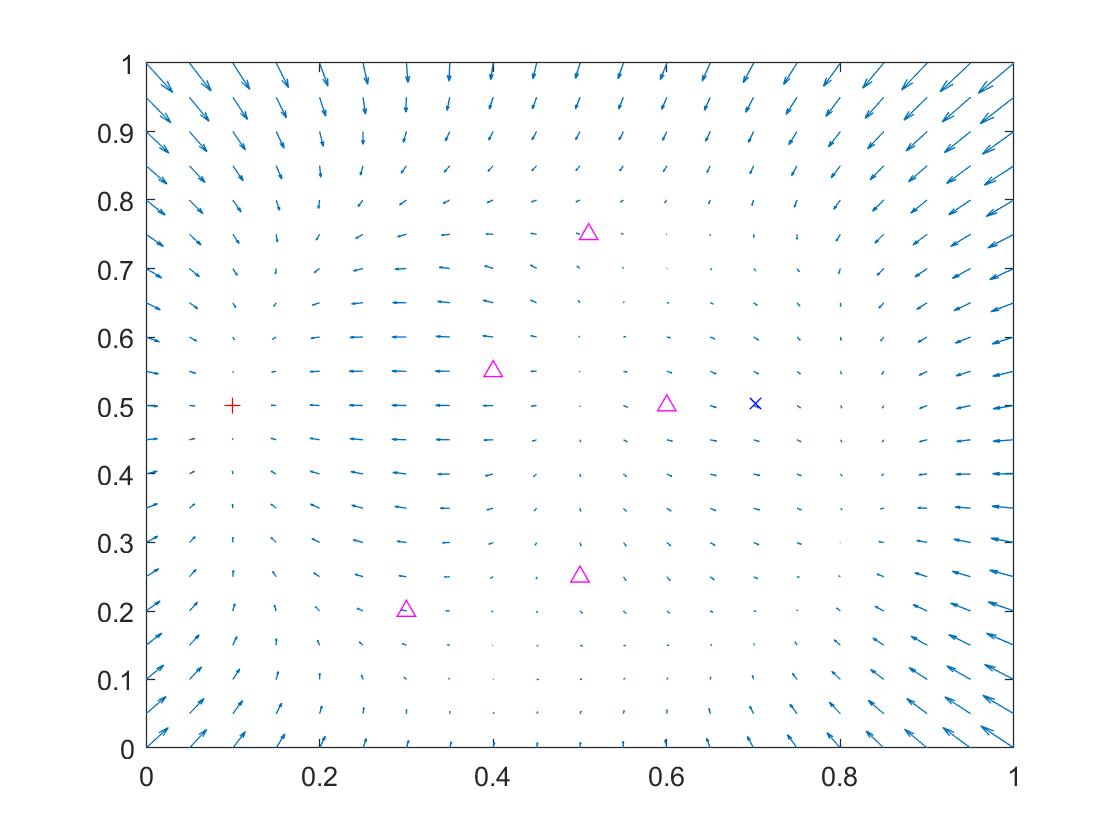}
    \caption{Vector field for $\nabla U(x)$. Sensor locations are given by magenta $\triangle$'s, target location (and global minimum of $U$) is given by red $+$, and mean initialization $\vx_{\text{avg}}(1)$ is given by blue $\times$.}
    \label{fig:vec-field}
  \end{subfigure}
  \begin{subfigure}[h]{0.45\textwidth}
    \includegraphics[width=\textwidth]{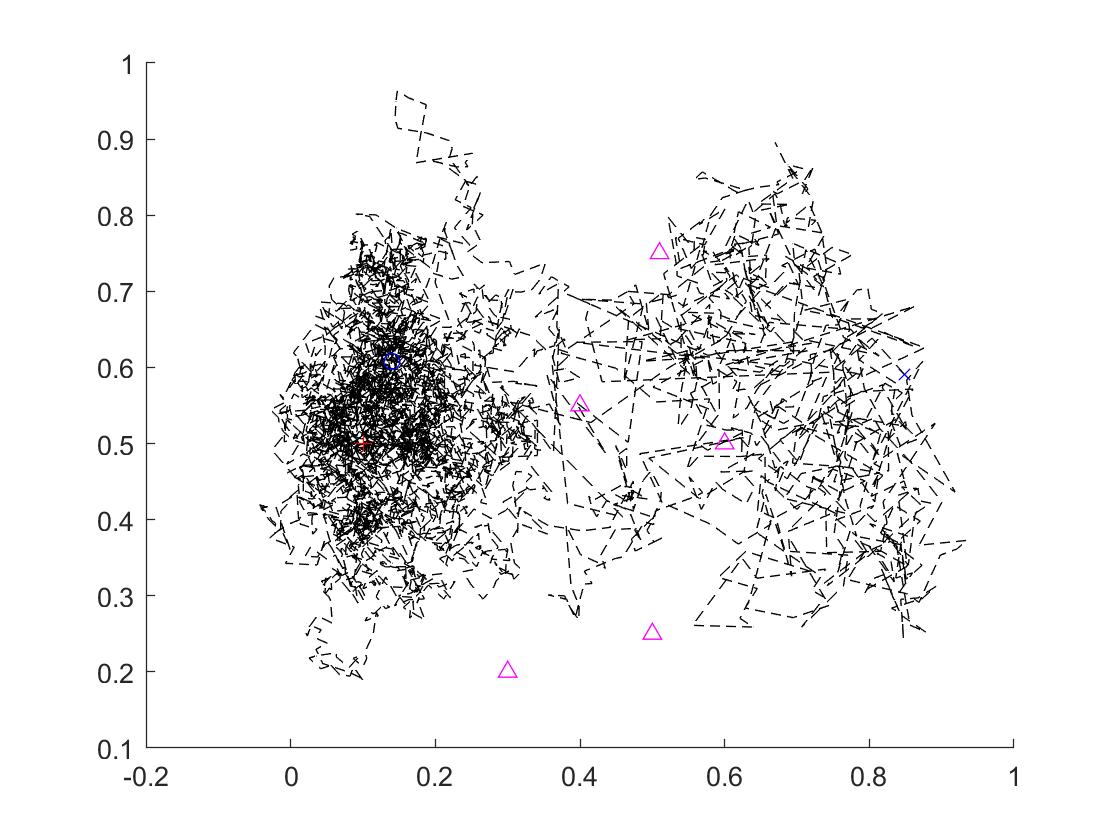}
    \caption{ Sample path of the distributed annealing algorithm. Magenta triangles denote sensor locations and red $+$ denotes target location. As $t\to\infty$, the process concentrates in the basin of the global minimum.}
    \label{fig:trajectory}
  \end{subfigure}
  \caption{}
\end{figure}

Note that, since we are dealing with only one target we have $d=2$ so that $U$ maps from $\R^2$ to $\R$. We leverage this low dimensionality to aid in visualizing the action of the algorithm. The gradient vector field $\nabla U(x)$ is plotted in Figure~\ref{fig:vec-field} along with the sensor and target locations. We emphasize that the vector field displayed in Figure~\ref{fig:vec-field} is the gradient vector field for $\nabla U(x)$ and not $\nabla \left(\sum_{n=1}^N U_n(x_n) \right)$. However, it is useful in visualizing the action of the algorithm since, as $\vx(t)$ approaches the consensus subspace, the average process $\vx_{\text{avg}}(t)$ asymptotically follows this vector field.\footnote{More precisely, the average process $\vx_{\text{avg}}(t)$ may be seen as an Euler discretization of the differential equation $\dot \vx = -\nabla U(\vx) + \vr(t),$ where $\vr(t)\to 0$ as $t\to\infty$.}

The unique global minimum of $U$ occurs at $x=z$, where $z$ is the target location. The vector field has a local minimum occurring near the point $(.8,.3)$ and multiple small-gradient regions that hamper the functioning of traditional gradient descent techniques.

We ran 100 trials of the algorithm for $10^4$ iterations each using the following weight parameters: $\alpha_t = 40\frac{1}{t}$, $\beta_t = .3\frac{1}{t^{1/4}}$, $\gamma_t = \frac{1}{(t\log(\log(t)))^{1/2}}$. To focus on the effects of annealing noise alone, we set $\bzeta_t \equiv 0$. Each trial used the same initial condition. The results of the simulations are displayed in Table~\ref{table:table1} and Figure~\ref{fig:vec-field}. Table~\ref{table:table1} considers the distance of the average $\vx_{\text{avg}}(t) = \frac{1}{N}\sum_{n=1}^N \vx_n(t)$ from the target location at various time instances. The table shows the number of trials for which $\vx_{\text{avg}}(t)$ fell within the ball $B_r(t)$ at (precisely) the iteration $t$ indicated in the column header. Theorem~\ref{th:main_result} implies that, for any $r>0$, the probability that $\vx_{\text{avg}}(t)$ lies inside the ball $B_r(z)$ about the target goes to 1 as $t\to\infty$. This is reflected in Table~\ref{table:table1}. We note that while we have not attempted to optimize convergence rate here, this may be a useful direction for future research.

\begin{table}[h]
\centering
\begin{tabular}{|c|c|c|c|c|c|c|}
\hline
 & $t = 500$ & $t=10^3$ & $t=2\times10^3$ & $t=5\times 10^3$ & $t=10^4$ \\ \hline
$r=.05$ & 8  & 10  & 13  & 14 & 18 \\ \hline
$r=.1$  & 29 & 26 & 39 & 41 & 50 \\ \hline
$r=.15$ & 44 & 45 & 52 & 56 & 72 \\ \hline
$r=.2$  & 59 & 62 & 70 & 71 & 84 \\ \hline
$r=.25$  & 69 & 69 & 75 & 83 & 89 \\ \hline
\end{tabular}
\caption{ The number of trials for which $x_{\text{avg}}(t)$ was in the ball of radius $r$ about the target (row) at iteration $t$ (column).}
\label{table:table1}
\end{table}

Figure~\ref{fig:trajectory} shows an example of a sample trajectory for a single trial after $5\times 10^3$ iterations. The trajectory diffuses through the state space, over time concentrating in the basin of the global minimum.

\section{Conclusions} \label{sec:conclusion}
We considered an annealing-based algorithm for computing global optima in distributed nonconvex optimization problems. The convergence result for the algorithm relies on several technical assumptions. Simple techniques for verifying that the technical assumptions hold were presented alongside a distributed target localization example.

\bibliographystyle{IEEEtran}
\bibliography{dist_glob_opt}

\end{document}